# DISCUSSION: CONDITIONAL GROWTH CHARTS

By Mary Lou Thompson[1]

*University of Washington*

I will use the terms "reference centiles" or "centile charts," as the setting that I consider here is more general than that of "growth charts."

Longitudinal reference centiles over some measure of time (typically age) are almost always implemented repeatedly on the same individual. In this kind of setting the notion of conditional or adaptive centile charts is very appealing, particularly when the within-individual variability is much less than that between individuals. While marginal or unconditional centile charts are common in many areas of application, conditional charts are still rarely encountered and further methodological development in this area is to be welcomed. The flexibility of the quantile regression approach of Wei and He (WH), for instance in allowing the dependence on past history to vary across centiles, is most attractive, as are the rigor and scope of their consideration of the problem.

I do, nevertheless, want to make a few cautionary remarks. The first relates to regression quantiles in particular, the second concerns a constraint common to all existing methods of constructing conditional percentiles, and the third and final point addresses the use of centile charts for screening. To concretize the discussion, the following setting will be considered throughout: the measurement of interest is assumed to be diastolic blood pressure in pregnant women, monitored between weeks 16 and 36 of pregnancy. There is typically an initial dip in blood pressure over this period, followed by a rise toward the end of pregnancy.

**1. Bias and precision.** My experience with the use of marginal regression quantiles has been that they are readily and robustly fitted, with far less of the "fine-tuning" that is needed for distributionally based centile estimation. Nevertheless, the flexibility of quantile regression estimates may come at a

Received December 2005.

[1]Supported in part by National Institutes of Health–National Institute of Child Health and Human Development Grant R01-HD-32562.

cost—should an appropriate distribution be identified, distributionally based estimates may well be more precise.

To evaluate bias and precision in marginal and conditional centile estimates, a simulation study was carried out on a presumed cohort of 1000 pregnant women, where it was assumed that the women were scheduled to attend an antenatal clinic once in each of five pregnancy intervals, namely during the weeks of gestation ("gestational age"): [16, 20), [20, 24), [24, 28), [28, 32), [32, 36). The visit times for each woman were assumed to be independently uniformly distributed within each interval. It was further assumed that the marginal distribution of the diastolic blood pressure of the $i$th woman at gestational age $t$, $Y_{it}$, was lognormally distributed with parameters (of the underlying normal distribution): $\mu_t = 4.247 - 0.019(t/10)^2 + 0.006(t/10)^3$ and $\sigma_t = 0.1$, where the units of blood pressure are mmHg. The $j$th measurement on the $i$th woman, $Y_{it_j}$, conditional on the measurement in the previous interval, $Y_{it_{j-1}}$, was again assumed to be lognormal with conditional parameters $\mu_{t_j|t_{j-1}} = \mu_{t_j} + \rho(\ln(Y_{it_{j-1}}) - \mu_{t_{j-1}})$ and $\sigma_{t_j|t_{j-1}} = 0.1(1-\rho^2)^{0.5}$.

A first-order autoregressive model [AR(1)] was assumed across intervals, with $\rho = 0.6$. It was further assumed that the probability of a woman attending a clinic in each of the prescribed intervals was 0.8, so that, on average, 20% of measurements are missing in each interval and overall. This approximates the situation that one might observe in practice. Figure 1 shows the longitudinal median blood pressure under this model as well as a simulated longitudinal sample with true percentiles superimposed.

Marginal and conditional centile estimates were obtained for 500 such simulated cohorts, using both the quantile regression approach (QR) described by WH and the LMS procedure [1]. Because the logged blood pressure measurements are multivariate normally distributed (MVN), this is also an ideal setting for the maximum likelihood approach suggested by Thompson and Fatti [4]. This last approach requires that it be possible to transform the longitudinal path, conditional on time points and covariates, to multivariate normality. Cubic splines were used to model the intercept term for each of the regression quantiles in the QR approach, to model L, M and S, and to model the mean of the multivariate normal distribution. The cross-sectional variance of the log transformed blood pressure measurements is constant and the AR(1) correlation, $\rho$, was also modeled as a constant in the MVN approach. Stata 9.1 was used for all analyses.

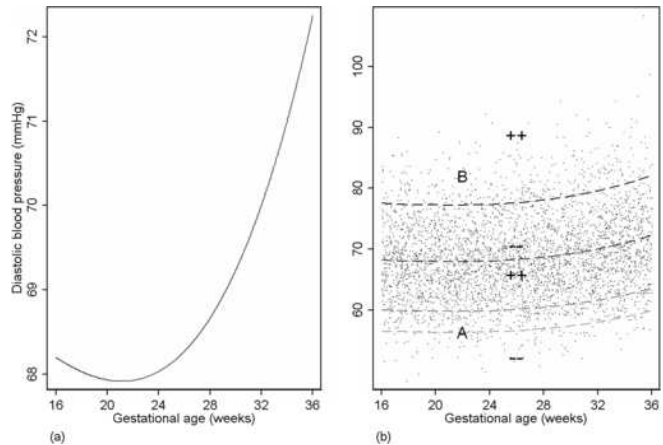

Fig. 1. (a) *Median and* (b) 3*rd,* 10*th,* 50*th,* 90*th,* 97*th percentiles.*

Note that simulations (not reported here) were also carried out using different numbers of subjects and measurement intervals and varying correlation, $\rho$. Results were consistent with those presented below.

All unconditional centile estimates were unbiased, but the variability in the quantile regression estimates was greater than that of the other two approaches, particularly at extreme percentiles; see Table 1.

Conditional centile estimates were considered at gestational age 26 for two hypothetical women, each with a previous measurement at week 22. The diastolic blood pressure reading in week 22 was assumed to be on the 3rd (marginal) percentile for the first woman, "A," that is, a blood pressure reading of 56.3 mmHg. The blood pressure of woman "B" in week 22 was assumed to be at the 97th percentile for that week, 82.0 mmHg. The measurements in week 22 and the true conditional 3rd and 97th percentiles for each woman in week 26 (A: "- -"; B: "++") are also shown in Figure 1.

The QR conditional centiles were estimated by fitting the model

$$Y_{i,j}(\tau) = g_{0,\tau}(t_{i,j}) + (\beta_{0,\tau} + \beta_{1,\tau}(t_{i,j} - t_{i,j-1}))Y_{i,j-1}$$

TABLE 1
*Standard deviation of estimates of marginal percentiles*

| Week | Method | 3rd | 10th | 50th | 90th | 97th |
|------|--------|------|------|------|------|------|
| 20   | QR     | 0.49 | 0.37 | 0.31 | 0.47 | 0.71 |
|      | LMS    | 0.40 | 0.30 | 0.28 | 0.39 | 0.57 |
|      | MVN    | 0.25 | 0.23 | 0.23 | 0.29 | 0.35 |
| 24   | QR     | 0.43 | 0.33 | 0.27 | 0.43 | 0.64 |
|      | LMS    | 0.36 | 0.28 | 0.25 | 0.37 | 0.53 |
|      | MVN    | 0.24 | 0.23 | 0.23 | 0.29 | 0.34 |
| 28   | QR     | 0.41 | 0.33 | 0.28 | 0.42 | 0.62 |
|      | LMS    | 0.35 | 0.28 | 0.25 | 0.36 | 0.52 |
|      | MVN    | 0.24 | 0.23 | 0.23 | 0.29 | 0.34 |
| 32   | QR     | 0.50 | 0.38 | 0.32 | 0.50 | 0.78 |
|      | LMS    | 0.42 | 0.32 | 0.28 | 0.41 | 0.62 |
|      | MVN    | 0.26 | 0.24 | 0.24 | 0.30 | 0.36 |

for $\tau = 0.03, 0.10, 0.50, 0.90, 0.97$, where $g_{0,\tau}(t)$ is modeled as a linear combination of five cubic basis splines. For the LMS approach, the longitudinal parameters L, M and S were estimated and then each observation was transformed to its corresponding $z$-score. Conditional centile estimates were then based on an AR(1) model, fitted to all $z$-scores that were one (visit) interval apart. The estimates of MVN conditional centiles were obtained by back-transforming to the observed scale the centile estimates from the Gaussian maximum likelihood estimates (see, e.g., [4]).

Table 2 summarizes the results of the simulation study comparing conditional centile estimates at week 26 under each of these approaches.

It can be seen that, except at the high (low) extreme percentiles for the low (high) prior path, the variability of the conditional QR and LMS estimates is similar here, but both are, not surprisingly, estimated with less precision than for the MVN-based model. The QR estimates are also slightly biased, as were QR estimates based on two previous observations (not shown). The form of the model for conditional centiles proposed by WH involves a linear adjustment for past history, whereas the lognormal conditional percentiles are extremely nonlinear in their relationship with previous measurements. It may well be that the extensions of the basic WH model, as discussed in Section 7.1 of their paper, would overcome this bias.

While regression quantiles do indeed provide an accessible and flexible means of estimating marginal and conditional percentiles, the above examples illustrate that gains can be made in terms of precision of estimates if an appropriate distributional structure can be identified. In addition, if bias is to be avoided, conditional percentile estimates using quantile regression will require careful choice of the form of the model for past history. As a counterpoint, however, as noted by WH, the conditional distributional structure may be more challenging to correctly identify than the marginal structure. WH provide an example of a setting where the distributional assumptions are not met and where hence the distributionally based centile estimates are biased.

TABLE 2
*Mean and standard deviation (SD) of simulated conditional percentile estimates at gestational age 26 weeks*

| Path "A" | 3rd (52.5)[a] | | 10th (55.1) | | 50th (61.0) | | 90th (67.6) | | 97th (70.9) | |
|---|---|---|---|---|---|---|---|---|---|---|
| **Method** | **Mean** | **SD** | **Mean** | **SD** | **Mean** | **SD** | **Mean** | **SD** | **Mean** | **SD** |
| QR | 52.8 | 0.57 | 55.5 | 0.42 | 61.5 | 0.33 | 68.3 | 0.51 | 71.8 | 0.79 |
| LMS | 52.5 | 0.55 | 55.1 | 0.43 | 61.0 | 0.33 | 67.6 | 0.37 | 70.9 | 0.41 |
| MVN | 52.5 | 0.19 | 55.1 | 0.20 | 61.0 | 0.22 | 67.6 | 0.27 | 70.9 | 0.31 |
| **Path "B"** | 3rd (65.8) | | 10th (69.0) | | 50th (76.4) | | 90th (84.7) | | 97th (88.8) | |
| QR | 65.3 | 0.61 | 68.7 | 0.46 | 76.2 | 0.37 | 84.5 | 0.57 | 88.6 | 0.91 |
| LMS | 65.8 | 0.38 | 69.0 | 0.40 | 76.4 | 0.44 | 84.7 | 0.67 | 88.8 | 0.91 |
| MVN | 65.8 | 0.29 | 69.0 | 0.29 | 76.4 | 0.28 | 84.7 | 0.30 | 88.9 | 0.32 |

[a]True conditional percentile (mmHg).

**2. Drift.** All of the methods of calculating conditional centiles implemented above can be expected to indicate once-off jumps in the path of an individual, but they are not able to deal satisfactorily with *drift* in an individual's path. WH acknowledge this when they state, for instance, that "The conditional growth charts may be unsuccessful in screening out subjects with gradual but persistent slowdown in growth."

To illustrate the problem in the pregnancy context, consider a woman, "C," whose blood pressure reading is on the 60th, 70th, 80th and 90th marginal percentiles in weeks 18, 22, 26 and 30, respectively. The discussion here will focus on true conditional percentiles, but the same comments carry over to the associated estimates under all three approaches considered here. In terms of the (true) conditional percentiles, the path "C" lies on the 68th, 74th and 83rd conditional percentiles in weeks 22, 26 and 30. The conditional centiles drift upward with the observations and thereby the woman's path is indicated as being progressively less extreme, relative to the marginal centiles.

Another plausible scenario would be one where there is a jump, after which the path remains steady at the new level. This might also be an indication of a potential problem, but unless the jump were large enough to be picked up immediately, the subsequent conditional percentiles would simply accommodate the change. Consider a woman, "D," whose blood pressure path has been moving along the (marginal) 50th percentile through week 22 of pregnancy and then at week 26 jumps to the (marginal) 80th percentile, where it remains in weeks 30 and 34. The reading in week 26 lies on the 85th conditional percentile but the subsequent readings in weeks 30 and 34 (both also on the 80th marginal percentile) lie on the 66th conditional percentile, because the conditional centiles have adjusted to the higher path.

In both examples the same features would be observed for the estimated conditional percentiles considered here. It is inherent in conditioning on past history that *all* past history is assumed "normal." One can argue that marginal and conditional centiles should be used together: for instance, the fourth observation on the above hypothetical woman "C" might be flagged at the 90th percentile of the marginal chart in week 30. However, these examples do illustrate a severe limitation in the usefulness of conditional centiles of this sort. Both of these types of paths are feasible in many contexts and these sorts of anomalies will not instill confidence in the nonstatistical user of such charts.

**3. Centile charts as a screening tool.** While centile charts may be of scientific interest in their own right, for example, to characterize or compare populations, they are generally constructed with a view to some sort of screening. In the context of children's growth considered by WH, they state, for instance, that: "When a measurement is extreme on the chart, the subject is often identified for further investigation. An extreme measurement is likely to be a reflection of some unusual underlying physical condition." If the intention is to use the conditional centile chart as a screen for a "problem" outcome such as future obesity in children's growth, or pre-eclampsia when monitoring blood pressure in pregnancy, then, as with any other screening method, it should be assessed in terms of its diagnostic accuracy.

There are several aspects of screening accuracy which merit attention here. First, note that there can be a strong association between a variable (e.g., a child's weight) and outcome (e.g., obesity) without the variable necessarily being a useful screening tool (see, e.g., [2]). Some examples drawn again from the blood pressure in pregnancy setting may illustrate this point.

Assume that a blood pressure reading at gestational age 22 weeks is the basis for conditional percentiles at week 26, which are being used to screen for some problem outcome ("disease"). Assume further that the distribution of diastolic blood pressure in the "diseased" is the same as that in the "normal" group prior to week 26, but that the two groups deviate in their means (but not variances) at week 26, at which point the percentage difference between the means for the "diseased" and "normal" groups is $d$. Then it is easily seen that, for a given specificity, $x$ (where specificity corresponds to the centile that is being used as a screen), the sensitivity is given by

$$\text{Sensitivity}(x) = \Phi\left(\frac{\ln(1+d)}{\sigma\sqrt{1-\rho^2}} - \Phi^{-1}(x)\right),$$

where $\Phi$ is the standard normal distribution function. In the above example, with $\rho = 0.6$, for there to be 90% sensitivity and specificity at week 26, there would have to be a 23% difference in mean blood pressure in the "diseased" relative to the "normal" population in that week. This represents an absolute difference in means of 15.6 mmHg, corresponding to a substantial 2.3 standard deviation difference in means. Note also that, if there was indeed no separation in the "diseased" and "normal" populations prior to week 26, then the screening sensitivity in those earlier weeks would be 1-specificity.

If, on the other hand, the percentage difference between the mean blood pressure paths of "diseased" and "normal" were a constant, $d$, at all gestational ages (all other distributional characteristics remaining the same), then it is also easily seen that

$$\text{Sensitivity}(x) = \Phi\left(\frac{\ln(1+d)\sqrt{1-\rho}}{\sigma\sqrt{1+\rho}} - \Phi^{-1}(x)\right).$$

Here, for a given specificity, sensitivity *decreases* with increasing $\rho$ and the marginal percentiles would have greater sensitivity than the conditional percentiles at all specificities. In the example described here, with $\rho = 0.6$, for there to be sensitivity and specificity of 90% in any given week, the percentage difference in means would need to be 67%. This example reinforces the discussion in Section 2, that conditional centiles are perhaps most useful in identifying jumps in a path.

The above examples consider once-off use of conditional centiles as a screen for "disease." In practice, of course, conditional centiles may be calculated at several points in time. With repeated screenings, overall sensitivity will increase at the expense of overall specificity (see, e.g., [3]). It should also be noted that, if the prevalence of the condition being screened for is low, the majority of "screened positive" individuals will be false positives. Depending on the consequences of a positive screen, this may have a range of sequelae ranging from emotional trauma to unnecessary invasive procedures.

WH have made a valuable contribution to the methodology available for estimating conditional centiles. However, there are limitations to the usefulness of such centiles as a screening mechanism. The implementation of centile charts, conditional or marginal, needs to be viewed in its entirety and this should include an evaluation of their screening effectiveness. It may well be that there is diligent measuring being carried out, to no useful purpose.

**Acknowledgment.** The assistance and insight of May M. Boggess, Stata Technical Services, with the Stata programming toward this discussion was much appreciated.

DEPARTMENT OF BIOSTATISTICS
UNIVERSITY OF WASHINGTON
SEATTLE, WASHINGTON 98195
USA
E-MAIL: mlt@u.washington.edu